\documentclass[11pt]{article}
\usepackage[greek,english]{babel}
\usepackage[round]{natbib}
\usepackage[continuous]{pagenote}
\makepagenote

\renewcommand*{\pagenotesubhead}[1]{}

\defcitealias{PR}{PR}
\defcitealias{RFM}{RFM}	
\defcitealias{LFM}{LFM}	
\defcitealias{NEM}{NEM}	

\makeatletter
\date{}
\makeatother

\begin{document}

\title{Wittgenstein, Peirce, and paradoxes of mathematical proof}

\author{Sergiy Koshkin\\
 Department of Mathematics and Statistics\\
 University of Houston-Downtown\\
 One Main Street\\
 Houston, TX 77002\\
 e-mail: koshkins@uhd.edu}
\maketitle
\begin{abstract} Wittgenstein's paradoxical theses that unproved propositions are meaningless, proofs form new concepts and rules, and contradictions are of limited concern, led to a variety of interpretations, most of them centered on the rule-following skepticism. We argue that his intuitions rather reflect resistance to treating meaning as fixed content, and are better understood in the light of C.S. Peirce's distinction between corollarial and theorematic proofs. We show how Peirce's insight that ``all necessary reasoning is diagrammatic", vindicated in modern epistemic logic and semantic information theory, helps explain the paradoxical ability of deduction to generate new knowledge and meaning. 
\bigskip

\bigskip


\end{abstract}

\section*{Introduction}

In his middle and late periods Wittgenstein reached conclusions that sound paradoxical, and at variance with mathematical practice as commonly understood. In Philosophical Remarks he argued that unproved propositions are meaningless and there can not be two different proofs of the same proposition, and in Remarks on the Foundations of Mathematics -- that proofs form new concepts, and that axioms and inference rules do not determine the theorems. 

Dummett framed Wittgenstein's position in terms of the {\it puzzle of deduction}, the tension between its two features, ``that in virtue of which we want to say that it yields nothing new, and that in virtue of which we want to say the opposite" \citep[p.\,299]{Dum73}\pagenote{On Dummett's reading, Wittgenstein's position takes on a Heraclitean or Hegelian flavor. According to  \citep[p.\,312]{Pap}, ``the Hegelian logic is not a solution of [Zeno's] paradox but a dismissal of the logical coordinates that generate it". Compare to Dummett's: ``Holism is not, in this sense, a theory of meaning: it is the denial that a theory of meaning is possible." \citep[p.\,309]{Dum73}.}. Other  interpretations of Wittgenstein's paradoxes have also been offered, some of them are briefly reviewed in Section \ref{SecInf}. What we would like to do is not to offer a yet another interpretation of Wittgenstein, but rather to understand the mathematical phenomena he highlights from a novel perspective. 

Wittgenstein and his interpreters largely treated all proofs as being of a kind, but C.S. Peirce, who pondered the puzzle of deduction a century earlier, distinguished between corollarial (routine) and theorematic (creative) proofs. 
In the 1960-s Hintikka rediscovered some of Peirce's ideas on mathematical proofs in his epistemic (modal) logic, and used them to resolve the puzzle of deduction. More recently, broadly Peircian approach to meaning and interpretation of mathematical proofs has been developed in the semantic information theory \cite{D'Ag}. We will argue that, despite the dissimilarities between the two thinkers\pagenote{The relationship between Peirce's and late Wittgenstein's positions is complicated. ``Meaning is use" is reminiscent of the pragmatic maxim (but qualified as ``sometimes, but not always"), and ``a way of grasping a rule that is not an interpretation" is akin to Peirce's habit change analysis. However, a detailed examination of the available evidence in \citep[Ch.1]{Bon} concludes that ``Wittgenstein expresses a basically negative attitude towards pragmatism as a {\it Weltanschauung}, but acknowledges affinities with pragmatism as a method". It is known that Wittgenstein read James extensively, and spent a year (1929) working with Ramsey, who developed his own version of semantic pragmatism based on Peirce's early works \citep{Mar}. Boncompagni speculates that Wittgenstein read Peirce's collection Chance, Love, and Logic, the Ramsey's source, some time after 1929. Ramsey was also a precursor of epistemic logic, with key ideas developed around 1929.}, Peirce's view of concepts and conceptual change in mathematics fits Wittgenstein's intuitions better than conventionalist, intuitionist or dialetheist interpretations, and largely defuses the charges of ``radical conventionalism" and ``assault on pure mathematics". But it also reveals some flaws in his analysis.

The first two sections discuss the first paradox in the inferentialist framework, characteristic of the Wittgenstein's middle period, and its role in his later abandonment of inferentialism. We turn to the second paradox and its diverse interpretations in Section \ref{Sec2Ap}. In Section \ref{SecCorTh} Peirce's corollarial/theorematic distinction is introduced, and related to the modern discussions of informal proofs and the informativity of deduction. In Section \ref{SecWittTh} we use it to argue that many of the Wittgenstein's theses are independent of the rule-following skepticism, and can be construed as rejection of the traditional idealization of conceptual omniscience, found also in Peirce's philosophy. In Section \ref{SecThPar}, motivated by Levy's refinement of the Peirce's distinction, we turn to a class of proofs that we call paradigmatic, which manifest conceptual shifts most explicitly. A model of mathematics, inspired by the modern epistemic logic, is sketched in Section \ref{SecEpiMath}, and it fits well with Wittgenstein's and Peirce's views on contradictions, reviewed in the following section. We summarize our discussion in Conclusions.

\section{No two proofs of one proposition}

The first paradox originates in the Wittgenstein's middle period, when he already believed that the  meaning of mathematical propositions is determined by their use, but interpreted this use as use in a ``calculus" \citep[p.\,201]{Rod97}. To avoid confusion, we will call a codified system for doing calculations and/or deductions a {\it formalism}. Wittgenstein's reasoning can be reconstructed as follows:

\begin{itemize}

\item[\bf P)] Meaning is use, and use in a formalism is use for inferring.

\item[\bf Q1)] Proposition is meaningful if it is inferentially linked to the axioms (proved), or if there is a decision procedure for producing such linkage\pagenote{There is some oscillation on Wittgenstein's part, noted in \citep[p.\,99]{Pleb}, as to whether merely having a decision procedure is enough to give meaning.}.

\item[\bf Q2)] Unproven propositions without a decision procedure are meaningless, and ``there can not be two independent proofs of one mathematical proposition" (\citetalias{PR}\pagenote{Standard abbreviations are used for Wittgenstein's works: PR for Philosophical Remarks, RFM for Remarks on Foundations of Mathematics, and LFM for Lectures on Foundations of Mathematics.}, 1975, p.\,184).

\end{itemize}

\noindent The reasoning from  {\bf Q1} to {\bf Q2} is as follows. A proof alters a formalism by turning a string of symbols into a usable proposition, it is the proof, or its blueprint at least, that enables its use and makes it meaningful. Another proof of the ``same" proposition will alter the meaning yet further, will link the sentence to different groups of axioms and/or in different ways, hence the proposition proved will not be the same. It is only our habit of attaching ``shadowy entities", meanings, to all well-formed sentences, even those that do not have any use, that leads us to believe in the sameness.

The fact that his conclusion is at odds with the common sense, and the common use of language, came to be unwelcome in the late period of ``philosophy leaves everything as it is". Late Wittgenstein replaced the ``calculi" as meaning givers by language games, and the rule-following considerations involved in them made the previously transparent notion of inference in a formalism problematic. But this by itself does not counter the logic of the no-two-proofs argument, if anything it makes it even stronger. Proofs are no longer rigid inferential chains, but performances, whose utility relies on reproducibility of the rule-following. But unproved propositions are still unusable, and hence meaningless. 

And yet in a remark from 1939-40 we read:``Of course it would be nonsense to say that {\it one} propistion can not have two proofs -- for we do say just that" (\citetalias{RFM}\pagenote{There are two different editions of RFM cited in the literature, with different numbering of the remarks. We cite the MIT paperback edition, as does Wright, but not Rodych and Steiner.}, II.58). Wittgenstein still seems to be torn between his old conception and the emerging late outlook, for he adds, ``proof is a mathematical entity that can not be replaced by any other; one can say that it can convince us of something that nothing else can, and this can be given expression by us assigning to it a proposition that we do not assign to any other proof" (\citetalias{RFM}, II.59). And in II.61 comes the crucial question:``How far does the application of a mathematical proposition depend on what is allowed to count as a proof of it and what is not?"

Wittgenstein's own answer comes in remarks from 1941:
{\small \begin{quote} It all depends on {\it what} settles the sense of a proposition, what we choose to say settles its sense. The use of the signs must settle it; but what do we count as the use? - That these proofs prove the same proposition means, e.g.: both demonstrate it as a suitable instrument for the same purpose. And the purpose is an allusion to something extra-mathematical (\citetalias{RFM}, V.7).
\end{quote} }
\noindent This singles out a sense of a proposition that remains unaltered throughout the play of linkages involved in different proofs, namely the sense bestowed by extra-mathematical applications. Hence, ``concepts which occur in `necessary' propositions must also occur and have a meaning in non-necessary ones" (\citetalias{RFM}, V.41). 

\citep[p.\,3]{Stein} argues that towards the end of 1930-s Wittgenstein's thought underwent a ``silent revolution", where he came to see mathematical propositions as ``hardened" empirical regularities, empirical generalizations a l\'a Mill promoted to the dignity of ``inexorable" rules. The same idea was expressed earlier in \citep[p.105]{Wri}, and it seems to be amply supported by multiple passages in RFM and LFM. The ``hardening" explains the stable reproducibility of the rule-following, and the widespread agreement on the outcomes of calculations and deductions, as well as applicability of formalisms to empirical matters, from which they were hardened. 

According to a number of scholars\pagenote{See e.g. \citep[p.329]{Moore}, \citep[p.\,28]{Pleb}, \citep[p.\,218]{Rod97} and 
\citep[p.\,23]{Stein}.}, this new stance had a bonus, perhaps part of the motivation for adopting it, of grounding Wittgenstein's hostility to mathematical logic and the upper reaches of set theory. During the middle period he could only fault them, or rather their (mis)interpretations, for assimilating extravagant formal games under the familiar concepts like numbers, sets, etc. Now he could say more, as in the oft-quoted \citetalias{RFM}, IV.2:
{\small \begin{quote}I want to say: it is essential to mathematics that its signs are also employed in mufti. It is the use outside mathematics, and so the meaning of the signs, that makes the sign-game into mathematics. Just as it is not logical inference either, for me to make a change from one formation to another... if these arrangements have not a linguistic function apart from this transformation.
\end{quote} }
\noindent Moore finds this passage to be ``essentially an assault on the very idea of pure mathematics" \citep[p.329]{Moore}. Thus, it seems that Wittgenstein's own solution left him at even greater conflict with mathematical practice than the no-two-proofs aporia it was meant to resolve.  

\section{Inferentialist solution}\label{SecInf}

Wittgenstein took the unwelcome conclusion of the first paradox as a strike against its premise, the equating of meaning to use in a formalism, in other words, he took his argument to be unsound. But, as reconstructed at least, it is invalid. Even if we identify meaning with inferential use, there is a problem with passing from {\bf P} to {\bf Q1}. Sure enough, the traditional rebuttal that comes to mind begs the question against Wittgenstein. We would like to say that we understand an unproved sentence by understanding its constituent parts and how they are linked. This appeals to content theories of meaning and the compositionality of language they support. But for Wittgenstein this, at best, transplants what applies to the empirical segment of language onto the grammatical sentences of mathematics, exactly the conceptual confusion he combatted in his middle and late periods. Dummet also objected that ``if Wittgenstein were right... communication would be in constant danger of simply breaking down" \citep[p.339]{Dum59}. But inferentialists do offer accounts of how languages can be mastered non-compositionally \citep[p.336]{Brand}, and communication rarely turns on nuances of meaning, as the utility of dictionaries indicates.

Still, there is no need to leave inferentialism behind to make sense of unproved sentences. There is even no need to compose them from simpler pieces occurring in other propositions, whose proofs are already known. If one wishes to use an unproved sentence inferentially one can assume it as a premise, and see what can be inferred from it. This is what Saccheri and Lambert did with the negation of the parallel postulate, and it gave them some idea of its meaning (enough for Saccheri to remark that it is ``repugnant to the nature of straight lines"). And this is what mathematicians continue to do with  odd perfect numbers or the Riemann hypothesis. Conversely, one can look for other unproved sentences, from which the one in question can be deduced, or better yet, for ones deductively equivalent to it. This was Sierpi\'nski's project for the continuum hypothesis. Of course, all such results are conjectural, but they do show that the inferential role does not reduce to a proof from the axioms. Moreover, if and when a proof or disproof of the sentence is found these conjectural results will be converted into proven or disproven propositions, and their proofs will quite literally contain the conjectural inferential chains as parts. Thus, the meaning of the proposition will ``contain" the meanings known before the proof even on the inferentialist conception.

This is not to say that ``the meaning" stays the same before and after the proof. By the same argument, its inferential role grows considerably. First, the proof establishes new inferential connections among different sentences of the formalism, and second, it delivers a new tool for proving other propositions. But the latter does not require one to even be familiar with the proof, just knowing (trusting) that there is a proof is enough. It is a common practice among mathematicians to make use of results they do not know proofs of. To summarize, quite a bit of use can be made of a proposition in a formalism independently of its proof. So, on the middle Wittgenstein's own terms, his argument is flawed.

But then whatever support it gave to altering its premise is also gone. Of course, Wittgenstein was no longer an inferentialist, so he may have had independent reasons for insisting on extra-mathematical use. One such reason is hinted at in RFM IV.25,``understanding a mathematical proposition is not guaranteed by its verbal form... The logical notation suppresses the structure". To Wittgenstein, the ``disastrous invasion of mathematics by logic", that masks conceptual leaps under the deceptive cover of familiar verbiage, is a target persisting through changes from the Tractatus to RFM. 

As for the extra-mathematical use, that put off Moore and many others, Wittgenstein is, in fact, quite equivocal. After the mufti quote, he goes on to ask: ``If the intended application of mathematics is essential, how about parts of mathematics whose application -- or at least what mathematicians take for their application -- is quite fantastic?...  Now, isn't one doing mathematics none the less?" (\citetalias{RFM}, IV.5). The answer comes in RFM, V.26 from a year later, and it is not what one might expect:
{\small \begin{quote}I have asked myself: if mathematics has a purely fanciful application, isn't it still mathematics? -- But the question arises: don't we call it `mathematics' only because e.g. there are transitions, bridges from the fanciful to non-fanciful applications?... But in that case isn't it incorrect to say: the essential thing about mathematics is that it forms concepts? -- For mathematics is after all an anthropological phenomenon. Thus we can recognize it as the essential thing about a great part of mathematics (of what is called `mathematics') and yet say that it plays no part in other regions... Mathematics is, then, a family; but that is not to say that we shall not mind what is incorporated into it.
\end{quote} }
This is hardly ``an assault on pure mathematics". In fact, it is reminiscent of Quine's division of mathematics into applied, its ``rounding out", and ``recreational". It seems that for Wittgenstein the use in mufti is just a check on the ``prose" surrounding the higher logic and set theory. But as we saw, neither the first paradox nor conceptual concerns make such use necessary.

\section{Proofs as rule-makers}\label{Sec2Ap}

In the late period Wittgenstein shifts to a much more diffused view of meaning than inferential use in a formalism. Accordingly, proofs are taken to grow a pre-existing meaning rather than to create it {\it ex nihilo}, and their contribution is framed in a broader context of language games. This leads to the second paradox. 

\begin{itemize}

\item[\bf P)] Proofs form new concepts and lay down new rules.

\item[\bf Q1)] In a proof we ``win through to a decision", placing it ``in a system of decisions" (\citetalias{RFM}, II.27).

\item[\bf Q2)] Formalism does not determine its theorems.

\end{itemize}
The main work is clearly done by the premise, and Wittgenstein amasses considerable amount of evidence to support it in RFM and LFM, see \citep[pp.\,39-40]{Wri} for a review. However, there is little consensus as to interpreting this premise, 
because, on the traditional views, it appears to be plainly false. In his influential 1959 interpretation Dummett denounced it as ``radical conventionalism":
{\small \begin{quote}He appears to hold that it is up to us to decide to regard any statement we happen to pick on as holding necessarily, if we choose to do so. [...] That one has the right simply to {\it lay down} that the assertion of a statement  of a given form is to be regarded as always justified, without regard to the use that has already been given to the words contained in the statement, seems to me mistaken \citep[p.337]{Dum59}.
\end{quote} }
\noindent On the Dummett's reading, Wittgenstein is even more radical than Quine, for whom holding on to a statement ``come what may" at least involves ``adjustments elsewhere in the system". But, as \citep{Str} pointed out, according to Wittgenstein, most mathematicians are usually compelled to accept a theorem when presented with a proof. This can hardly be compared to laying down a convention.

Wright remarks that ``it ought to be possible, after we have accepted the proof, satisfactorily to convey what our understanding of a statement used to be", and concludes that it is not, in fact, possible on the traditional accounts of meaning as content. Because if a proof conforms to the old content it can not also create a new one \citep[pp.\,53-54]{Wri}. He then suggests that Wittgenstein's talk of ``conceptual change" is figurative, and is meant to dislodge the traditional figure of ``recognizing" what our rules already dictate, which is generally the target of the rule-following considerations. Wittgenstein's figure comes with figures of speech, like ``inventions" and ``decisions" in place of ``discoveries" and ``recognitions", and is meant to play a therapeutic role (Ibid. pp.\,48-49).

However, as we saw with the first paradox, it is possible, {\it pace} Wright, to give an account of a meaning before and after the proof, which makes sense of meaning change without appealing to the rule-following. It involved giving up the view of meaning as content, even intuitionist content. One can make sense of the change even on content theories, but such a change will be, as Dummett put it in his modified ``more plausible" reading of 1973\pagenote{Dummett reaffirmed and elaborated on his modified reading in \citep{Dum94}, which reproduces some passages from his 1973 lecture almost verbatim.}, banal. A new characterization of an ellipse, say, would give us a new rule for recognizing that something is an ellipse, which we did not have before the proof. But ``the new criterion will always agree with the old criteria, when these are correctly applied in accordance with our original standards... even if we failed to notice the fact" \citep[p.53]{Dum94}. He then suggests that a robust interpretation of Wittgenstein's thesis requires an example in which the old and the new criteria disagree, while we are unable to find any mistakes, either in the proof or in the application of the criteria, a seemingly impossible feat. We would have to claim that the mistake is there even if we are unable, in principle, to locate it. Only an all-seeing God can then distinguish the banal and the robust interpretations, and rejecting such an Olympian view is exactly the Wittgenstein's point, according to Dummett.

\citep{Stein} gives a yet another interpretation, somewhat reminiscent of Stroud's, based on the view that he attributes to Fogelin. On this view, we observe widespread agreement on what constitutes following a rule, because the rules themselves are empirical regularities promoted to the dignity of a rule, ``hardened", ``{\it because} they all agree in what they do we lay it down as a rule and put it down in the archives" (\citetalias[XI]{LFM}). This Copernican turn throws a new light on the before and after of a proof. Professionals, trained as they are in the ways of their language game's rule-following, will be particularly compelled to accept a proved proposition as the only possible outcome. But this itself is an empirical regularity of behavior after training. And empirical regularities do break down, training is not destiny. Hence, what proof delivers, while not a legislated convention, falls short of a foregone conclusion. Before the proof, Wittgenstein continues in LFM: 
{\small \begin{quote}The road is not yet actually built. You could if you wished assume it isn't so. You would get into an awful mess. [...] If we adopt the idea that you could continue either in this way or in that way (Goldbach's theorem true or not true) -- then a hunch that it will be proved true  is a hunch that people will find it the only way of proceeding.
\end{quote} }
This should give some idea of the diversity of opinion on the issue, but note that most of it revolves around the role of the rule-following. After looking deeper into the puzzle of deduction we will see that the rule-following may not be the only issue.

\section{Corollarial/theorematic distinction}\label{SecCorTh}

Peirce's self-described ``first real discovery about mathematical procedure" was a generalization to all deductive reasoning of a traditional distinction between the ``logical" and ``geometric" consequences in Euclidean geometry, traceable as far back as Aristotle. The former can be read off of the diagram directly, while the latter require auxiliary constructions, ``which are not at all required or suggested by any previous proposition, and which the conclusion... says nothing about" (\citetalias{NEM}\pagenote{NEM v:p is a standard abbreviation for The New Elements of Mathematics by Charles S. Peirce, v volume, p page.}, 4:49). Earlier, the distinction inspired Kant's distinction between the analytic and synthetic arguments. Most of Peirce's writings on the subject remained unpublished until 1970-s, so the distinction remained buried until Hintikka brought it back from obscurity in 1979 \citep[p.56]{Hin80}, after rediscovering a version of it in his own work. 

Peirce developed a diagrammatic version of the first order predicate calculus with quantifiers (existential graphs), which allowed him to argue that ``all necessary reasoning is diagrammatic" and extended the corollarial/theorematic distinction to all deductions \citep[p.56]{Dip}.
{\small \begin{quote}Of course, a diagram is required to comprehend any assertion. My two genera of Deductions are 1st those in which any Diagram of a state of things in which the premisses are true represents the conclusion to be true and such reasoning I call corollarial because all the corrollaries that different editors have added to Euclid's Elements are of this nature. 2nd kind. To the diagram of the truth of the Premisses something else has to be added, which is usually a mere May-be and then the conclusion appears. I call this theorematic reasoning because all the most important theorems are of this nature. (\citetalias{NEM}, 3:869)
\end{quote} }
\noindent Peirce's view is supported by the modern studies of diagrammatic reasoning \citep[p.24ff]{Giaq}. But after Frege, along with the construction generally, the distinction came to be seen as ``psychologistic", and in geometry specifically, as an artifact of its incomplete formalization. 

Peirce characterizes a theorematic proof as introducing a ``foreign idea, using it, and finally deducing a conclusion from which it is eliminated" (\citetalias{NEM}, 4:42). This foreign idea is ``something not implied in the conceptions so far gained, which neither the definition of the object of research nor anything yet known about could of themselves suggest, although they give room for it" (\citetalias{NEM}, 4:49). Theorematic reasoning reflects the informal idea of mathematicians about the non-triviality of proofs. In contrast, corollarial reasoning is routine, and is closely related to what middle Wittgenstein called a ``decision procedure". However, even in theories with effective (algorithmic) proof procedures the actual proving of theorems may not be routine, because the procedures are too complex, and, therefore, intractable. For example, the elementary Euclidean geometry and Boolean algebra are effectively decidable, but their general decision procedures are intractably complex. The algorithmic complexity of deductions correlates with their informativity \citep[p.175]{D'Ag}. If one thinks of information as, in Hintikka's slogan, elimination of uncertainty, then one can see how theorematic proofs are informative\pagenote{Hintikka's extensional view of uncertainty was rather narrow compared to Peirce's own. A general critique of his interpretation of the corollarial/theorematic distinction is \citep{Ket85}.}. They eliminate genuine uncertainty about what they prove (Ibid., p.178), whereas corollarial (tractably algorithmic) proofs do not. 

Several measures of informativity/complexity of (formal) deductions have been proposed in the modern epistemic logic and semantic information theory. The first one was Hintikka's {\it depth}, the number of new layers of quantifiers introduced in the course of the proof. It is also motivated by the auxiliary constructions in geometry, Hintikka analogizes them to the new ``individuals" introduced when the newly quantified variables are instantiated in natural deduction systems. However, Hintikka's depth does not detect all types of theorematic steps. They appear even in proving Boolean tautologies, where no quantifiers are present, but extra letters and/or connectives are introduced in the intermediate formulae \citep[p.62]{Dip}. In response, D'Agostino and Floridi proposed to supplement it with a second depth, which is in play even in proving Boolean tautologies. It is the depth of nested patterns of subarguments that introduce and discharge additional assumptions in a natural deduction system \citep[p.178]{D'Ag}. Jago proposed a single alternative measure, the shortest proof length in a sequent calculus without the contractions and the cut \citep[p.331]{Jag13}. As these explications show, informativity is relative to the background formalism, and incremental -- the depths or proof lengths depend on a chosen proof system, and mark heap-like changes rather than sharp divides.

Even so, formalization and measurement of qualitative change can only go so far. Informativity within a formal system invites the picture where the conceptual resources are circumscribed in advance, and deductions simply spread truth values to some previously undecided propositions. This is a picture adopted by Dummett. In the same lecture where he modified his interpretation of Wittgenstein, he insists on what we will call {\it conceptual omniscience}. It is a semantic version of Hintikka's logical (better to say, epistemic) omniscience, the idealization that the knowledge of premises entails the knowledge of all of their deductive consequences. Proofs do grow knowledge, according to Dummett, but not meaning. That they can not do while staying faithful to the prior content of propositions, they merely facilitate verification of other claims, mathematical or empirical. Deduction brings new knowledge 
{\it despite} preserving the meanings \pagenote{Dummett's solution to the puzzle of deduction is criticized in \citep{Haack}.}. And this is enough to affirm a strong form of deductive determinism: once the axioms are laid down the theorems are determined for everyone but a radical skeptic about the rule-following. 

However, as already Peirce pointed out, theorematic reasoning involves ``foreign ideas", concept formation or transformation over and above the theorem's formulation, and the background knowledge. The nature of these new concepts is suggested by his examples, and is made explicit in the modern semantic information theory. They manifest in the construction and/or recognition of new patterns, auxiliary figures in geometry, composite structures in set theory, or compound predicates and propositional formulae in formal systems \citep[p.170]{D'Ag}. One defines new objects, and/or finds new ways to describe their properties and interrelations with other objects, old and new. Many proved properties are turned into new definitions. The conceptual omniscience is problematic because much of mathematicians' effort goes into {\it crafting} definitions, and few theorems are proved about objects introduced already in the axioms. The skeletal semantics of the model theory, that parses formulae down to the basic elements, is not the semantics of informal proofs \citep[p.18]{Azz}. To use Dummett's own example, the concept of ellipse does not appear in either planimetric or stereometric axioms, and it is only one among an infinite variety of objects they give room for. That theorems about ellipses should be proved at all is not determined by the formalism. 

Of course, ellipses are strongly motivated by common observations, but this suggests exactly the empirically mediated ``determinacy" that Wittgenstein describes. In the practice of mathematics, definitions do more than single out formal patterns. The newly formed concepts are linked to concepts from other formalisms, informal intuitions, and applications outside of mathematics. When the conceptual resources are specified in advance, the interpretational labor required to make proofs and theorems meaningful can not be captured by them. And ``without an interpretation of the language of the
formal system the end-formula of the derivation says nothing; and so nothing is proved" \citep[p.26]{Giaq}. The meaning of unproved theorems is not determined because, after all, we may not be {\it smart enough} to deduce them, let alone anticipate the concepts to be introduced in their proofs, or statements. The appearance of elliptic curves and modular forms in the Wiles's proof of the Last Fermat theorem gives an idea of just how much concept formation can be involved. 

While the informativity detects (in degrees) the need for concept formation, it does not express it. The Peirce's theorematicity is intended to capture the accompanying conceptual surplus, which emerges even when working in completely formalized deductive systems. Thus, {\it pace} Dummett, we can make a non-banal sense of how proofs form new concepts and rules without offering the impossible counterexamples to proved theorems.

\section{Wittgenstein and theorematic proofs}\label{SecWittTh}

As we argued, proofs can effect conceptual change even aside from the rule-following indeterminacy. The irony is that not only the commentators tended to overlook the corollarial/theorematic distinction, but so did late Wittgenstein himself. The difference is that if they, in effect, treated all deductive reasoning as corollarial, he treated it all as theorematic. Middle Wittgenstein admitted, at least occasionally, that effective decision procedures give sense even to unproved propositions:``We may only put a question in mathematics (or make a conjecture) where the answer runs:``I must work it out"" (\citetalias{PR}, p.151). But late Wittgenstein dropped the distinction in favor of a uniform approach. This approach might have, indeed, caused a radical breakdown in communication, between him and his interpreters. If a proof effects conceptual change no matter what kind of proof it is, one needs a conception of this that applies to all cases, and one might as well analyze the simplest cases, corollarial ones. 

While most of Wittgenstein's examples are theorematic\pagenote{Examples of proofs discussed in RFM include: conversion of strokes into decimals, occurrence of 770/777 in the decimal expansion of $\pi$, impossibility of listing fractions in the order of magnitude, impossibility of angle trisection with straightedge and compass, recursive abbreviations in Principia, Cantor's diagonal argument, identification of real numbers with Dedekind cuts, and G\"odel's incompleteness theorem.}, he is also fond of stressing the equivalence between a formalism and a calculus, deduction and calculation. On Peirce's view, the essential difference is that calculation (as in adding and multiplying numbers) involves no theorematic steps, one just works it out. But, at the same time, the distinction is relative and incremental, so Wittgenstein might have seen no philosophical ground to draw a sharp line in the sand. 

Whatever his reasons, Wittgenstein forced his interpreters to fit his conceptual change thesis even to the most routine of calculations, and to explain how to conceive of it when the informativity of deduction all but disappears. And this invariably left the general rule-following indeterminacy as the only viable option, see e.g. \citep[pp.\,48--49, 145--147]{Wri}. In hindsight, one can see how applying even the paper-and-pencil addition algorithm to nevertofore seen numbers has a residue of theorematicity to it. Because who is to say that the addition as previously grasped is not really quaddition, and so 68+57=5 \citep[p.\,9]{KrW}. But without the benefit of examples where the ``foreign idea" is more substantive, it is easy to miss the non-banal residue. The addition algorithm has been mechanized since the first arithmometers, and one needs thick skeptical glasses to discern conceptual change in adding 68 to 57. Wright arrives at something like this infinitesimal theorematic residue reading, when drawing the contrast to the more substantive case of the Last Fermat Theorem:
{\small \begin{quote}All that doing number theory does is acquaint us with a variety of constructions  which are deemed analogous... A proof of Fermat's theorem, if we get one\pagenote{Wright was writing in 1980. Wiles first announced his proof in 1993, but it contained a gap. The final version, completed in collaboration with Taylor, did not appear until 1995.}, may not closely mimic these other constructions; it may rather appeal to a general concept which they illustrate, and then present new methods as relevant to it... In contrast, we can circumscribe the technique relevant to the solution of some problem of effectively decidable type absolutely exactly \citep[p.55]{Wri}.
\end{quote} }
Late Wittgenstein might have (legitimately) taken exception to the ``absolutely", but, perhaps, it would have better served his ends to offer a sop to Cerberus\pagenote{In a 1908 letter to lady Welby Peirce explains his description of a sign as having effect upon a person as follows:``My insertion of ``upon a person" is a sop to Cerberus, because I despair of making my own broader conception understood".}, instead of ignoring the contrast altogether. As it is, even Wright only gives the above interpretation in the context of ascribing to Wittgenstein the intuitionist semantics of proofs (p.\,54), and later uses the same intuitionist gloss in discussing the occurrence of 777 in the decimal expansion of $\pi$ (p.145). There he remarks that the amount of uncertainty about such occurrence ``contrasts with the scope which we should expect occasionally to have for discretion" (p.150), if we only had loose analogies.

Thus, we are left with the general rule-following skepticism directly applied to the decimal expansion of $\pi$. But such skepticism infects any discourse, including empirical assertions that Wittgenstein pains to distinguish from mathematical ones. If ``the further expansion of an irrational number is a further expansion of mathematics" 
(\citetalias{RFM}, IV.9) means that genuine discretion can be exercised in deciding whether 777 occurs or not, Wittgenstein is in trouble. But, as the texts quoted by Fogelin and Steiner suggest, this is not what it means. The absence of a reason for the rule-following is not a reason for the absence of the rule-following. Wittgenstein did not deny that the rule-following in proofs typically produces a determinate result, he argued that the traditional accounts misconstrue the nature of this determinacy. 

In short, the corollarial/theorematic perspective explains away and/or accommodates the diverging interpretations of the second paradox, and dulls its edge in the process. Its conclusion is revealed to hold for all proofs only legalistically\pagenote{Commenting on his provocative early assertion that ``any statement can be held true come what may", Quine writes in Two Dogmas in Retrospect: ``This is true enough in a legalistic sort of way, but it diverts attention from what is more to the point: the varying degrees of proximity to observation...".}, substantively only for theorematic proofs, and even then not in the sense of leaving room for genuine discretion required by Dummett for non-banality. Still, this is only a part of the story.

\section{From theorematic to paradigmatic}\label{SecThPar}

As we saw, distinguishing corollarial and theorematic proofs helps contextualize Wittgenstein's theses, and move the focus away from the rule-following. But theorematic proofs are not created equal either. Levy pointed out that under the heading of theorematic reasoning Peirce describes a wide range of examples \citep[p.\,99]{Levy}. On one end, we have Euclid's auxiliary lines, and clever algebraic substitutions; on the other, Fermat's ``infinite descent" (mathematical induction), and Cantor's diagonal argument applied to general power sets. Theorematicity comes in degrees, but in the two latter cases the historical context suggests more than a difference in degree. Euclid and Cardano were applying already established axioms\pagenote{Of course, even in the case of Euclid, ``axiom" in the modern sense applies only loosely.} of geometry and algebra, while Fermat, and especially Cantor, were introducing new ones. 

Levy describes the distinction as the one between using ideas logically implied by the principles already adopted, perhaps tacitly, and ideas demanding the adoption of new principles (Ibid.). Let us follow him in splitting off proofs appealing to such new principles, which we will call {\it paradigmatic}, the word often used by Wittgenstein himself. In Peirce's terms, paradigmatic proofs appeal to something not only unimplied by conceptions so far gained, but to something they do not even give room for. This is an informal analog of the difference between conservative and non-conservative extensions of a formal theory. A conservative extension introduces new concepts and principles in such a way that their use can be eliminated from  proofs, as long as they are absent from the theorems' statements. In a non-conservative extension, previously undecidable propositions may become provable  \citep[p.\,20]{Azz}. For example, a strengthened form of Ramsey's theorem about graph colorings, due to Parris and Harrington, is undecidable in the first order arithmetic, but is provable in the ZFC set theory. 

Of course, what the principles give room for is somewhat open to interpretation, unless they are completely formalized. 
In informal practice, the theorematic/paradigmatic boundary is blurred, for mathematicians rarely work within a fixed formal system. The Wiles's, or Parris-Harrington's, proofs were not seen as paradigmatic (in the narrow sense), because modern number theorists do not confine their paradigm to the first order arithmetic. In these terms, theorematic proofs extend the theorem's background, albeit conservatively (in the broad sense), while corollarial ones do not. 

The theorematic/paradigmatic divide also parallels Toulmin's distinction between the warrant-using and warrant-establishing arguments in the argumentation theory, for which he invokes Ryle's metaphor of traveling along a railway already built versus building a new one \citep[p.\,120]{Toul}. He also points out that, historically, ``deductions" referred to all warrant-using arguments, not only to the formal logical ones. They included, for example, astronomers' calculations of eclipses based on Newton's theory, and Sherlock Holmes's surmises from crime scene evidence, which certainly involved theorematic reasoning. The parallel with Wittgenstein's own metaphors of building ``new roads for traffic" (\citetalias{RFM}, I.165), ``designing new paths for the layout of a garden" (\citetalias{RFM}, I.166), and ``building a road across the moors" (\citetalias{LFM}, X) should be plain. Except for Wittgenstein, {\it every} proof ushers in a new paradigm, he distinguishes paradigmatic from theorematic no more than theorematic from corollarial, at least not explicitly. 

Most of the commentary tacitly assumes that ``the proofs" are proofs in modern-style formalisms, with explicitly stated axioms and rules of inference. But most of Wittgenstein's examples in RFM involve historical proofs produced in no such formalisms. Moreover, in RFM, II.80 he explicitly states: ``It is often useful in order to help clarify a philosophical problem, to imagine the historical development, e.g. in mathematics, as quite different from what it actually was. If it had been different {\it no one would have had the idea} of saying what is actually said" [emphasis added]. Let us look at some of Wittgenstein's examples in this light.

That the angle trisection is possible by neusis (with {\it marked} straightedge and compass), was known in antiquity, and that Euclid would rule out such constructions was not determined by the loose idea of straightedge and compass. Similarly, identifying Dedekind cuts with the real numbers was not determined by special real numbers, and vague generalities about them, known before Dedekind. Indeed, the prevailing conception of the continuum was Aristotelian, on which it is not assembled from points/numbers at all. Wittgenstein charges that Dedekind established a new rule for what a real number is under the misleading cover of a familiar geometric cut:``The division of rational numbers into classes did not {\it originally} have any meaning, until we drew attention to a particular thing that could be so described. The concept is taken over from the everyday use of language and that is why it immediately looks as if it had to have a meaning for numbers too" (\citetalias{RFM}, IV.34). The cut is exactly a composite structure, a new pattern, generally implicated in the concept formation through proofs. Moreover, as we now know, even arithmetized continuum does not have to consist of Dedekind cuts, the real numbers, it could instead be hyperreal or the absolute continuum of Conway, both containing infinitesimals. 

The Cantor's diagonal argument brings in a controversial at the time idea of actual infinity, and an even more controversial idea of comparing such infinities according to the Hume's principle of bijective correspondence. Even Bolzano, Cantor's precursor, rejected the Hume's principle because it conflicted with the Euclid's part-whole axiom (the whole is greater than its part) for infinite sets \citep[p.625]{Manc}. G\"odel gave an influential argument that Cantor-style cardinalities were inevitable as measures of infinite size, but alternative measures that preserve the part-whole axiom, so-called numerosities, were later found nonetheless \citep[p.637]{Manc}. Wittgenstein surmises that instead of emphasizing the disanalogy between the real and the natural numbers, that the diagonal argument brings out, the cardinality talk reduces it to mere difference in size. Again, ``the dangerous, deceptive thing" is ``making what is determination, formation, of a concept look like a fact of nature" (\citetalias{RFM}, App.\,II.3). 

Wittgenstein, it seems, has a case to resist the idealization of conceptual omniscience, whether he intended to make it or not. In the case of paradigmatic proofs, not only are the Dummett's impossible counterexamples not needed, they are, in fact, possible. One might object that only complete formalization fixes the meaning of concepts, and in paradigmatic cases we are dealing with informal proofs operating with loose concepts. But this is how mathematics evolved historically: we did not have formal concepts {\it prior} to a proof, and had it conform to them, formalisms were developed {\it after}, if not as a result of, the proof's adoption. The prior use involved concepts, such as they were, that were consistent with the adoption of conflicting alternatives. If one of them is then adopted, what is it if not a conceptual change? This would accord well with Peirce's habit-based view of the meaning of concepts:

{\small \begin{quote}The concept which is a logical interpretant is only imperfectly so. It somewhat partakes of the nature of a verbal definition, and is as inferior to the habit, and much in the same way, as a verbal definition is inferior to the real definition. The deliberately formed, self-analyzing habit -- self-analyzing because formed by the aid of analysis of the exercises that nourished it -- is the living definition, the veritable and final logical interpretant. Consequently, the most perfect account of a concept that words can convey will consist in a description of the habit which that concept is calculated to produce. (CP\pagenote{CP v.p is a standard abbreviation for The Collected Papers of Charles Sanders Peirce, v volume, p paragraph.}, 5.491)
\end{quote} }

Of course, Dummett saw Wittgenstein as talking about proofs in a modern formalism, and he might concede the change introduced by paradigmatic proofs as again a banal point, that is what makes them paradigmatic. Fair enough. But the determinacy is often claimed even for paradigmatic cases, as with the Cantor's cardinalities, and this claim is then relied upon to present proving in a formalism as a model, a cleaned up version, a ``rational reconstruction", as Carnap and Reichenbach called it, of how mathematical knowledge is acquired. Moreover, as we argue next, the paradigmatic shades into the theorematic just as the theorematic shades into the corollarial.

\section{Epistemic model of mathematics}\label{SecEpiMath}

What might an alternative model of mathematical development, more hospitable to Wittgenstein's intuitions, look like? It will be helpful to frame the changes induced by proofs in terms of epistemic logic. A formalized version of such a picture is developed in \citep[p.329]{Jag09}\pagenote{Jago conceives of the epistemic horizon very differently, and abstracts from the informal shell. In Conclusions, we explain why his formal framework may also be unattractive to Wittgenstein due to the conceptual omniscience concerns.}. 

At any given time only some propositions of the formalism are known (proved). Not even all of their corollarial consequences can be said to be known, not because there is a problem with deducing them, but because there may be no reason to turn attention to them. When an occasion arises, say in applications, they will be deduced as a matter of routine. We may even take some low grade theorematic reasoning (below a vaguely marked threshold) as part of the routine, this resembles what Kuhn called the ``normal science" of ``puzzle-solving". There are also propositions, like intermediate formulae in cumbersome computations, that are only significant in the context of deducing something else, and would not be attended to on their own. They may be corollarial, but even if they already occurred in known proofs they may not be portable enough to register as independent items of knowledge. They only become epistemically relevant when one is working through a known proof, or attempting a new one. 

What we have, then, is an epistemic core of theorems surrounded by a desert, epistemic horizon, of unclaimed and/or technical propositions, through which passage to any (truly) new theorem lies. The core is immersed into an informal shell of motivations, analogies, interpretations, and applications, that supplement the meaning of concepts featured in it, and may, occasionally, even conflict with the formalism. But whatever the formalism does express conceptually, is largely limited to its epistemic core. The shell motivates some anticipations and hunches extending beyond it, and some non-core parts may be explored -- by deriving antecedents and consequents of some conjectures, and exploring new concepts and techniques that show promise.

We can now better appreciate the similarities and the differences between theorematic and paradigmatic proofs. Both will expand the epistemic core and constrain the informal shell, by sorting conflicting intuitions and providing new rules for the ``puzzle-solving". A theorematic proof will do so conservatively, making the new rules seem like validations of prior commitments. A paradigmatic proof, in contrast, will have to {\it negotiate} the axioms already adopted, and the informal anticipations of the shell. This is how it was with the Cantor-Dedekind arithmetization of the continuum, or with the Zermelo's well-ordering proof. Of course, a theorematic proof may reveal that the formal terms conflict too much with their informal counterparts (as almost happened with Zermelo's proof). However, if anything is rejected in such a situation it will not be the proof itself, but rather the formalism, at least on the traditional account.

There is a problem with that account, however. We can {\it legislate} that accepting a proof always counts as ``conforming" to prior rules, and altering the formalism counts as ``modifying" them, but this convention is at odds with historical practice. A foreign idea in theorematic proofs may be treated as transgressing the rules, rather than as conforming to them, for the rules may not have been meant to be applied {\it this way}. Conversely, a proof may induce a conceptual shift even if it accords with the previously adopted rules. 

The Weierstrass's example of a continuous nowhere differentiable function caused a shift in understanding continuity, even though it followed the already adopted formal definition of Cauchy. Presumably, to conform to prior concepts one would have had to change the formalism. This illustrates how a formalism's ability to fix the concepts does not extend far beyond its epistemic core. Uninterpreted formal theorems may be {\it syntactically} determined by the formal transcription rules, but, as such, they are conceptually thin, ``understanding a mathematical proposition is not guaranteed by its verbal form". And conceptualized theorems are not fixed by the formalism alone, and therefore are not determined by it. The case for determinism turns not (merely) on the rule-following, but on the conceptual omniscience, without it the Wittgenstein's thesis is defensible.

The syntactic idealization is at odds even with the Platonist and intuitionist accounts, where the formalisms do not fully capture semantic consequence and mathematical truth, the very accounts that motivate content theories of meaning. The axiom of replacement was added to Zermelo's original axiomatization of set theory because the latter was seen as inadequate to express the Cantorian ``inductive conception of sets". The subsequent search for large and larger cardinals indicates that even ZFC does not fully capture that conception. In fact, {\it any} formalism, including the Euclidean geometry and Peano arithmetic, can not fully capture the ``intended" concepts on the Platonist or intuitionist interpretations of mathematics. Those belong to the platonic realm, or to the synthetic potential of a quasi-Kantian subject.

But if revision of formalisms need not amount to conceptual revision, then their affirmation need not amount to conceptual conformity either. And if so, every novel proof puts the formalism on the line and forces a decision one way or the other. Even if the proof is accepted, we still have a conceptual shift and a new rule, an extension of mathematics. Wittgenstein might have expressed himself thus: a formalism may determine {\it its} theorems (barring the rule-following indeterminacy), but not what they  mean, and a new proof may reveal that it failed to mean {\it our} concepts. Put this way, Wittgenstein's point is neither conventionalist nor banal, it is, indeed, a radical departure, but not from the mathematical practice. Rather, it is a departure from the prevailing philosophical prose of its rational reconstruction, which presupposes conceptual omniscience. 

\section{Ex falso nihil fit}\label{SecNih}

That late Wittgenstein's intuitions line up with the epistemic model of mathematics is further corroborated by his view of contradictions. If a formalism is inconsistent then, under the {\it ex falso quodlibet} rule, anything, literally, goes. But does this mean that an inconsistent formalism fails to capture any concepts? From the epistemic perspective, the only contradictions that affect practice are the known ones. Hidden contradictions, beyond the epistemic core, can not threaten the use of a formalism, and therefore do not preclude it from being conceptually meaningful. If a theorematic foreign idea leads to a contradiction we may take it as a sign that the formalism was no good, but we may also take it as a sign that the foreign idea was too foreign, and save (the consistent fragment of) the formalism by blocking its use. This is how Russell saved the Frege's system, by restricting the Basic Law V\pagenote{Basic Law V leads to the unrestricted comprehension and Russell's paradox.}. Wittgenstein's own example is arithmetic:``If a contradiction were now actually found in arithmetic -- that would only prove that an arithmetic with {\it such} a contradiction in it could render very good service; and it would be better for us to modify our concept of the certainty required, than to say that it really not yet have been  a proper arithmetic" (\citetalias{RFM}, V.28). And this explains his 
{\it ex falso nihil fit} proposal:``Well then, don't draw any conclusions from a contradiction. Make that a rule"(\citetalias{LFM}, XXI). 

While dialetheists do see Wittgenstein as a precursor \citep{PrR}, it does not seem that he had something like paraconsistent logic in mind. Paraconsistent logicians go to much trouble beyond the {\it ex falso nihil fit} to neutralize contradictions. This is because, as Turing already pointed out at one of Wittgenstein's lectures, any conclusions, derivable from a contradiction in a classical formalism, can also be derived without going through any contradictions. The rules of inference have to be altered quite dramatically to block all such derivations. 

This is only needed, however, if one insists on syntactic, mechanizable transcription rules. The Wittgenstein's ``rule" amounts instead to boxing the formalism within its prior epistemic horizon, where no contradictions arise. This consistent fragment stood, and was used, on its own, it is only the {\it ex post facto} projection of contradictions derived later that makes one think that there was anything wrong with it. ```Up to now a good angel has preserved us from going {\it this} way'. Well, what more do you want? One might say, I believe: a good angel will always be necessary whatever you do" (\citetalias{RFM}, II.81). In a way, this is Wittgenstein's dissolution of the Gettier problem of epistemic luck. 

A good angel, it is true, is already relied upon in assuming that training is effective and machines do not break down, but it still helps to take precautions. Reliability, like theorematicity, comes in degrees, and Wittgenstein is disregarding, it seems, the higher reliability of mechanizable rules, as opposed to an open-ended ``if I {\it see} a contradiction, then will be the time to do something about it" (Ibid.). What we do not use can not hurt us, he argues, and even when a contradiction comes to light -- ``what prevents us from sealing it off? That we do not know our way about in the calculus. Then {\it that} is the harm" (Ibid.). However, it is prudent to minimize the stumbling around even when we do not (yet) know our way about, and we know empirically that mechanizable rules are apt to accomplish that\pagenote{A telling example is the practice of the Italian school of algebraic geometry in Wittgenstein's lifetime under Enriques and Severi, who adopted a more laissez faire attitude to mathematical rigor, and relied on intuition to find their way about. The results produced by the Italians eventually became unreliable, and later had to be reworked in the formal framework of Weil and Zariski. Mumford wrote about Severi's 1935-1950 work: ``It is hard to untangle everywhere what he conjectured and what he proved and, unfortunately, some of his conclusions are incorrect" \citep[p.\,326]{Brig}.}. Therefore, they are preferable by the late Wittgenstein's own lights, it is only the  prose surrounding them that he can object to.

Ramsey, a presumed bridge between Peirce and Wittgenstein, anticipated some ideas of the epistemic ("human") logic in his papers written around 1929, when he worked with Wittgenstein at Cambridge. The passages on consistency quoted in \citep[p.71]{Mar} are quite suggestive:
{\small \begin{quote}We want our beliefs to be consistent not only with one another but also with the facts: nor is it even clear that consistency is always advantageous; it may well be better to be sometimes right than never right. Nor when we wish to be consistent are we always able to be: there are mathematical propositions whose truth or falsity cannot as yet be decided. [...] human logic or the logic of truth, which tells men how they should think, is not merely independent of but sometimes actually incompatible with formal logic.
\end{quote} }
Peirce's pragmatic attitude towards the hidden contradictions is also known \citep[p.237]{Murph}, it follows from his general rejection of the Cartesian ``paper" doubt. According to Peirce, mathematics generally has no need for formal logic, as its own method of ideal experimentation is more basic, and consistency of mathematical theories, like any other scientific claim, is to be doubted only when there comes up a specific reason to do so. And if it should happen, Peirce, like Wittgenstein, was confident that mathematicians will be up to the task of addressing it. However, Peirce was equally pragmatic about the usefulness of rigor and formal rules, indeed he developed a number of formal systems himself.

The tolerance of contradictions reinforces our earlier point about conceptual determinacy: if deducing a contradiction does not ``nullify" the original formalism the latter can not be said to determine its conceptual meaning simply by the syntactic consequence. Inconsistency is yet another symptom of the coming apart between formalisms and informal shells that make them meaningful.

\section{Conclusions}

We argued that the first paradox is aimed against the static theory of meaning, the semantics of fixed content. Unproved theorems are not quite meaningless, even on the inferentialist semantics, but their meaning grows with new proofs. Proving ``the same" proposition twice is like entering the Heraclitean river twice, -- it is not quite the same. The second paradox replaces inferentialism with a pragmatist, in spirit, semantics of rule-governed practice. That a formalism grounded in it determines its theorems can only be maintained if the formalism is assumed to have preconceived content, and to be  executed by clockwork subjects. Once these idealizations are dropped the indeterminacy of theorems loses the air of a paradox, even without the breakdown in the rule-following clockwork. The higher tolerance for contradictions also becomes more palatable in this de-idealized picture.

This is not to say that Wittgenstein's arguments are without flaws. Proofs bring conceptual change in degrees, noted already by Peirce, at the extremes of which we find mechanical corollarial proofs and trailblazing paradigmatic ones, with a theorematic continuum in between. While only the rule-following considerations make the corollarial conclusions indeterminate, the theorematic conclusions display genuine indeterminacy, due to the conceptual limitations of the formalism's users. Idealizing away these limitations, and the conceptual flux they create, leads to the puzzle of deduction's triviality, on the traditional accounts of mathematics. The semantics of preconceived content can only accommodate Wittgenstein's theses as banalities. Paraconsistent logic is still off the mark with its syntactic blocking of blatant contradictions that did not bother Wittgenstein. But perhaps hidden contradictions should have bothered him some more, in view of the pragmatic advantages that consistent formalisms provide when it comes to the ``use in mufti". 

Finally, while the traditional accounts overstate the conceptual determinacy and the fixity of meaning, the Wittgensteinian alternative faces the opposite problem. The puzzle of deduction remains, albeit turned on its head -- it is not the non-triviality of deduction that is puzzling, but rather its conformity \citep[p.\,301]{Dum73}. Wittgenstein's allusions to empirical regularities in this regard are intriguing, but obscure, relations between meaning, content and empirical regularities need further elaboration. There is a similar, and better understood, puzzle concerning the continuity of knowledge across Kuhn's scientific revolutions, which may provide some guidance.  

Epistemic logic offers an illuminating perspective on Wittgenstein's paradoxes, but the extensional turn it took in Hintikka's and subsequent work would likely make it unattractive to either Wittgenstein or Peirce. Jago, for instance, takes as a platitude Hintikka's thesis that epistemic growth amounts to ruling out possibilities, with the possibilities described in terms of (classically impossible) possible worlds \citep[p.329]{Jag09}. But this is only a platitude if one accepts that the possibilities are conceptually determined, and specifiable {\it in advance}. The development of knowledge can be described as narrowing down pre-existent options only if we are deploying concepts that will emerge before they actually do. This is exactly the conceptual omniscience, the Olympian view, that Wittgenstein took pains to oppose. 

Peirce offered an alternative approach, that Wittgenstein might have found more congenial. Instead of working with conceptually determined possibilities, like the possible worlds, he talked of constraints on them in terms of vague descriptions. Such constraints on future knowledge can be formulated even in terms of past concepts, without the Olympian view. How continuity of knowledge across scientific revolutions can be understood, along the Peircean lines, is sketched e.g. in \citep[p.\,274ff.]{Short}. A similar approach to mathematics seems promising. Unfortunately, intensional approaches to epistemic modality remain underdeveloped.

{\footnotesize
\printnotes
}

\nocite{CP}

{\small
\bibliographystyle{plainnat}
\bibliography{WittAporia}
}

\end{document}